\documentclass[10pt]{amsart}
\usepackage{amsmath,amsthm,amssymb,amscd}
\usepackage{epsfig}

\begin{document}
\newcommand{\m}{{\mathfrak m}}                             
\newcommand{\M}  [1] {\ensuremath{{\overline{\mathcal M}}{^{#1}_0(\R)}}}   
\newcommand{\cM} [1] {\ensuremath{{\mathcal M}_{0}^{#1}}}                  
\newcommand{\CM} [1] {\ensuremath{{\overline{\mathcal M}}{^{#1}_0}}}       
\newcommand{\oM} [1] {\ensuremath{{\mathcal M}_{0}^{#1}(\R)}}              
\newcommand{\dM} [1] {\ensuremath{{\widetilde{\mathcal M}}_{0}^{#1}(\R)}}  
\newcommand{\oQ} [1] {\ensuremath{{\mathcal Z}^{#1}}}                      
\newcommand{\Z}  [1] {\ensuremath{{\overline{\mathcal Z}}{^{#1}}}}         
\newcommand{\PGL} {\Pj\Gl_2(\R)}                           
\newcommand{\RP} {\R\Pj^1}                                 
\newcommand{\C} {{\mathbb C}}                              
\newcommand{\R} {{\mathbb R}}                              
\newcommand{\Int} {{\mathbb Z}}                            
\newcommand{\Pj} {{\mathbb P}}                             
\newcommand{\T} {{\mathbb T}}                              
\newcommand{\La} {{\mathbb L}}                             
\newcommand{\Sg} {\mathbb S}                               
\newcommand{\Gl} {{\rm Gl}}                                

\theoremstyle{plain}
\newtheorem{thm}{Theorem}[subsection]
\newtheorem{prop}[thm]{Proposition}
\newtheorem{cor}[thm]{Corollary}
\newtheorem{lem}[thm]{Lemma}
\newtheorem{conj}[thm]{Conjecture}

\theoremstyle{definition}
\newtheorem{defn}[thm]{Definition}
\newtheorem{exmp}[thm]{Example}

\theoremstyle{remark}
\newtheorem*{rem}{Remark}
\newtheorem*{hnote}{Historical Note}
\newtheorem*{nota}{Notation}
\newtheorem*{ack}{Acknowledgments}

\title {A Space of Cyclohedra}

\author{Satyan L. Devadoss}
\address{Department of Mathematics, Ohio State University, Columbus OH 43210}
\email{devadoss@math.ohio-state.edu}

\begin{abstract}
The real points of the Deligne-Mumford-Knudsen moduli space \CM{n} of marked points on the sphere has a natural tiling by associahedra.  We extend this idea to create a moduli space tiled by {\em cyclohedra}.  We explore the structure of this space, coming from blow-ups of hyperplane arrangements, as well as discuss possibilities of its role in knot theory and mathematical physics.
\end {abstract}

\maketitle

{\small
\begin{ack}
I am grateful to Mike Davis for pointing out the affine example given in Brown \cite{bro} and to Vic Reiner for conversations about non-crossing partitions.  Thanks also go to Jack Morava and Jim Stasheff for their encouragement.  Finally, I am indebted to the late Rodica Simion who provided motivation, discussions, and enthusiasm.  This paper is dedicated in her memory.
\end{ack}}


\baselineskip=15pt

%
%
\section{Introduction}

Introduced in algebraic geometry, the moduli space \cM{n} of Riemann spheres with $n$ labeled punctures has become a central object in mathematical physics.  There is a natural compactification \CM{n} of this space whose importance was emphasized by Grothendieck in his famous \emph{Esquisse d'un programme} \cite{gro}. It plays a crucial role in the theory of Gromov-Witten invariants and symplectic geometry, also appearing in the work of Kontsevich on quantum cohomology, and closely related to the operads of homotopy theory.  The {\em real} points \M{n} of this space, the set of points fixed under complex conjugation, will provide our inspiration.

A beautiful fact about \M{n} is its tiling by convex polytopes known as {\em associahedra} \cite{dev}.  This paper is motivated in part by Rodica Simion who asked the author whether there is an analogous space tiled by the {\em cyclohedron} polytope.  We show the existence of such a space, which we denote as \Z{n}, describing it in terms of compactifications of configuration spaces. One remarkable feature is that \Z{n} is a $K(\pi, 1)$ space (it is aspherical).

There are two approaches taken to construct this space:  One is local, coming from gluing copies of polytopes, giving a combinatorial feel to the problem.  The second is a global perspective, in terms of blow-ups of certain hyperplane arrangements.  The importance of \Z{n} comes from its cyclohedral tessellation.  This polytope first appears in the work of Bott and Taubes who introduce it in the context of non-perturbative link invariants \cite[\S1]{bt}.  It continues to show up in different areas of mathematical physics and especially knot theory~\cite{bar}.

%
%
\section {Tiling by Associahedra}

\subsection{}
Detailed ideas behind the motivating example of \M{n} and its tessellation by the associahedron are given in \cite{dev}.  But for the uninitiated, we present a quick review, highlighting the key points which will allow us to compare and contrast the cyclohedral space constructed in the next section.  We first give the classic construction.

\begin{defn}
The {\em associahedron} $K_n$ is a convex polytope of dimension $n-2$ with codimension $k$ faces corresponding to using $k$ sets of meaningful bracketings on $n$ variables.
\end{defn}

\begin{figure}[h]
\centering {\includegraphics {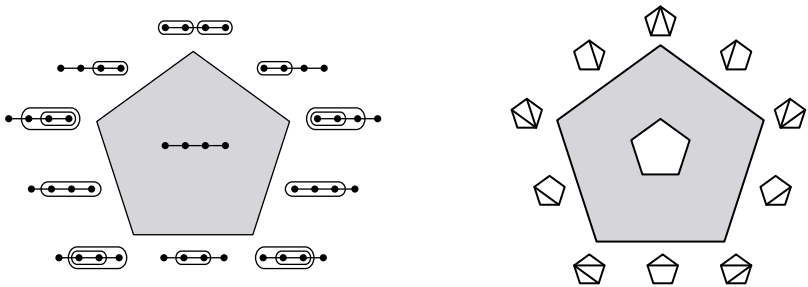}}
\caption{Associahedron $K_4$ using parentheses and polygons}
\label{k4pg}
\end{figure}

Stasheff originally defined the associahedron for use in homotopy theory in connection with associativity properties of $H$-spaces \cite[\S2]{jds}.  He was able to show sphericity, describing it as the face poset of a CW-ball.  It was given a realization as a convex polytope originally by Milnor (unpublished), and later by Lee \cite{lee} and Haiman (unpublished).  There is an alternate definition which we will base our work on, with Figure~\ref{k4pg} showing an example of the relationship:

\begin{defn}
The {\em associahedron} $K_n$ is a convex polytope of dim $n-2$ with codim $k$ faces corresponding to using $k$ sets of non-intersecting diagonals\footnote{Mention of {\em diagonals} will henceforth mean non-intersecting ones.} on an $(n+1)$-gon.
\end{defn}

\subsection{} \label{ss:compact}
We define the compactified moduli space \CM{n} of Riemann punctures. Although this variety is defined over the integers, we look at the real points:

\begin{defn}
The space \M{n} is the {\em Deligne-Knudsen-Mumford} compactification of the configurations of $n$ {\em labeled} points on $\RP$ quotiented by the action of $\PGL$.  That is,
$$\M{n} = \overline{{\rm Config}^n(\R \Pj^1)/\PGL}.$$
\end{defn}

\noindent Since $\PGL$ is three-dimensional, \M{n} is a manifold of dimension $n-3$. The moduli space is a point when $n=3$ and $\RP$ when $n=4$. For $n>4$, these spaces become non-orientable, getting extremely complicated as $n$ increases.

To understand how the associahedron tiles \M{n}, the idea of compactification needs to be explored.  Think of configurations of $n$ points on an $m$-manifold (in our case, this space is $\RP$), where points are allowed to move but not to touch.  Compactifying allows the points to collide and a {\em system} is introduced to record the {\em directions} points arrive at the collision.  In the celebrated work of Fulton and MacPherson~\cite{fm}, this method is brought to rigor in the algebro-geometric context.\footnote{Fulton-MacPherson is over the complex, but for our purposes the reals are more relevant.}  As points collide, they land on a {\em screen}, viewed as $\R\Pj^m$, which is identified with the point of collision.\footnote{These projective spheres have been dubbed {\em bubbles}, and the compactification process as {\em bubbling} (see \cite{dev2} for details).}  Now these points on the screen are themselves allowed to move and  collide, landing on higher level screens.  However, they are free up to an action of $\Pj\Gl_{m+1}(\R)$, the affine automorphism on each screen.  Kontsevich describes the process in terms of a magnifying glass:  On a particular level, only a configuration of points is noticeable; but one can zoom-in on a particular point and peer into its screen, seeing the space of the collided points.

\begin{figure}[h]
\centering {\includegraphics {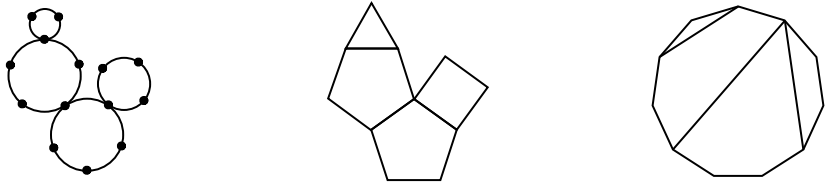}}
\caption{Bubbles and Polygons}
\label{bubblegon}
\end{figure}

The compactification of \oM{n} comes from Deligne-Knudsen-Mumford, closely related to the Fulton-MacPherson method.  There are roughly two distinctions:  The base space of configurations and the bubbles from collisions are both $\RP$ with $\PGL$ acting on them, making them indistinguishable.  As points collide, the result is a new bubble fused to the old at the point of collision, where the collided points are now on the new bubble.  Figure~\ref{bubblegon}a illustrates one such example; one cannot tell on which $\RP$ the bubbling began.
The second difference is that each bubble needs to have a minimum of three marked points.  The reason for this is due to the stability condition coming from Geometric Invariant Theory \cite[\S 8]{git}.  GIT gives a natural compactification of the space of genus zero algebraic curves which are {\em stable} in the sense of having only finitely many automorphisms; that is, the curves cannot have just two distinguished points.


\subsection{}
The duality presented in Figure~\ref{bubblegon} shows the relation between compactified space of points on the $\RP$ bubbles and polygons with non-intersecting diagonals.  Roughly, the $n$ points on a bubble correspond to an $n$-gon, with levels of collisions between the points giving rise to diagonals.  It naturally follows that $K_{n-1}$ appears in \M{n}, where the {\em tessellation} by copies of associahedra comes from dealing not just with points on bubbles, but {\em labeled} points.  In other words, by keeping track of possible labelings, we find \M{n} to be tiled by exactly $\frac{1}{2}(n-1)!$ copies of $K_{n-1}$:  Although there are $n!$ possible ways of labeling the sides of an $n$-gon, since the dihedral group $D_n$ (of order $2n$) is in $\Pj \Gl_2(\R)$, each associahedral domain of \M{n} corresponds, in some sense, to a labeled $n$-gon up to rotation and reflection.  The idea behind how these associahedra glue together becomes fundamentally combinatorial.

\begin{figure} [h]
\centering {\includegraphics {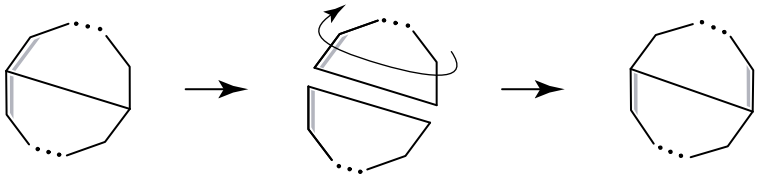}}
\caption{{\em Twist} along a diagonal}
\label{ktwist}
\end{figure}

\begin{defn}
The {\em twist} of a polygon $G$ along a diagonal $d$, denoted by $\nabla_d(G)$, is the polygon obtained by separating $G$ along $d$, `twisting' (reflecting) one of the pieces, and gluing them back (Figure~\ref{ktwist}).  It does not matter which piece is twisted since the two results are identified by an action of $\PGL$.
\end{defn}

\begin{thm} \textup{\cite[\S4]{dev}} \label{t:twist}
Two polygons $G_1, G_2$ with labeled sides, representing codim $k$ faces of associahedra, are identified in \,\M{n} if there exist diagonals $d_1, \ldots, d_r$\, of \,$G_1$\, such that $$(\nabla_{d_1} \cdots \nabla_{d_r}) (G_1) = G_2.$$
\end{thm}

\begin{figure} [h]
\centering {\includegraphics {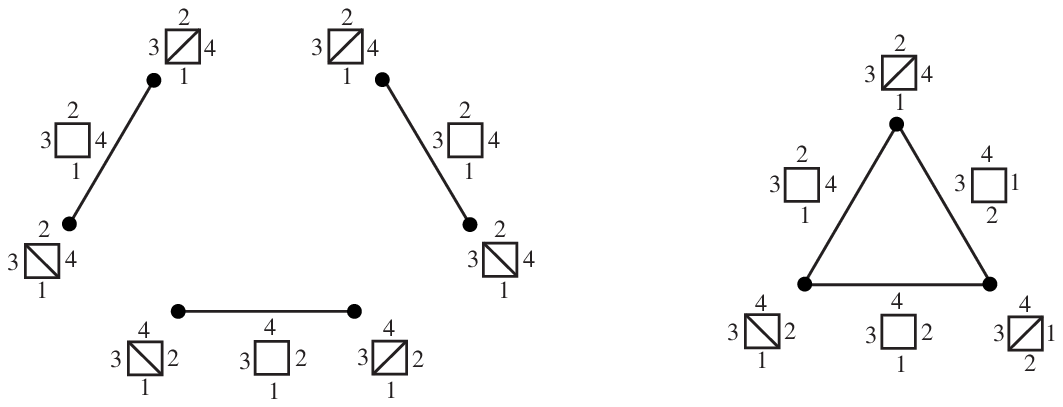}}
\caption{Copies of $K_3$ and \M{4}}
\label{m04pcs}
\end{figure}

Indeed, this shows how copies of associahedron glue together to form the moduli space. Since a codim $k$ face of $K_n$ has $k$ diagonals, then $2^k$ distinct associahedra meet at the face in the tessellation.  Figure~\ref{m04pcs} demonstrates \M{4} tiled by three labeled $K_3$'s. The {\em cross-ratio} is a homeomorphism from \M{4} to $\RP$, identifying three of the four  points with $0, 1,$ and $\infty$. The final point moves between the three fixed ones (interiors of $K_3$'s) or collides with them (endpoints where associahedra glue together).

%
%
\section{Tiling by Cyclohedra}

\subsection{}
The cyclohedron $W_n$ was originally manifested in the work of Bott and Taubes and later given its name by Stasheff.  The reason for a similar name is due to its construction which closely mimics $K_n$.  

\begin{defn}
The {\em cyclohedron} $W_n$ is a convex polytope of dim $n-1$ with codim $k$ faces corresponding to using $k$ sets of meaningful bracketings on $n$ variables arranged on a circle.
\end{defn}

As it is for $K_n$, it would be natural to expect an alternate definition in terms of dissected polygons.  It was the clever idea of Simion to come up with one.  She formulates it in terms of {\em centrally symmetric} polygons with even numbers of sides.  Henceforth, any mention of a (non-intersecting) diagonal on such a $2n$-gon will either mean  a {\em pair} of {\em centrally symmetric} diagonals or the diameter of the polygon.  

\begin{defn} \label{d:cyclo}
The {\em cyclohedron} $W_n$ is a convex polytope of dim $n-1$ with codim $k$ faces corresponding to using $k$ sets of non-intersecting diagonals on a \break $2n$-gon.
\end{defn}

\begin{figure}[h]
\centering {\includegraphics {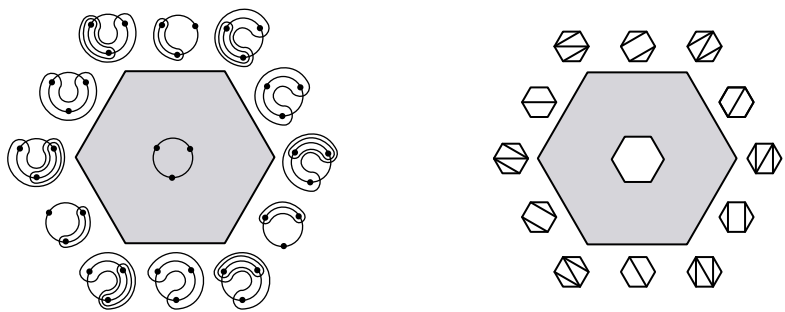}}
\caption{Cyclohedron $W_3$ using parentheses and polygons}
\label{w3pg}
\end{figure}

Figure~\ref{w3pg} shows an example of the two descriptions of $W_3$.  Note a crucial difference between $K_n$ and $W_n$:  For the codim one faces of the associahedron, we do not place brackets around all $n$ variables.  In contrast, the cyclohedron allows this since one can distinguish the cyclic manner in which the $n$ variables are combined.  This fact plays a role in the compactification defined below.


\subsection{}
Before moving to understand the space of cyclohedra, it is worthwhile to explore its properties.  It was noted by Stasheff \cite[\S2]{jds} that the boundaries of $K_n$ can be identified with products of lower dimensional associahedra.  He makes a similar observation for $W_n$.

\begin{prop} \textup{\cite[\S4]{jds2}}
Faces of cyclohedra are products of lower dimensional cyclohedra {\em and} associahedra.  In particular, given \,$W_n$ where \,$k-1+\sum n_i = n$, there are various inclusions of the form
$$W_k \times K_{n_1} \times \ldots K_{n_k} \hookrightarrow W_n.$$
\end{prop}

\begin{figure}[h]
\centering {\includegraphics {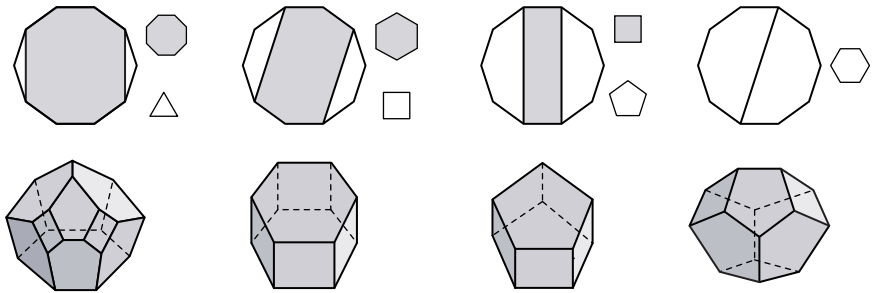}}
\caption{Codim one faces of $W_5$}
\label{w5codim1}
\end{figure}

\begin{exmp}
Since $W_5$ is a four-dimensional polytope, we look at its possible top dim (codim one) faces.  They are given by the different ways of placing a diagonal in a $10$-gon.  Figure~\ref{w5codim1} illustrates the four possible cases.  The far left is the product $W_4 \times K_2$.  Since $K_2$ is simply a point, the result is $W_4$ itself.  The middle figures show $W_3 \times K_3$ and $W_2 \times K_4$, where on one hand the line segment is $K_3$ and on the other $W_2$. The last possible type is $W_1 \times K_5$, which is simply $K_5$.
\end{exmp}

\begin{rem}
For the associahedron, it is well known that the inclusions of lower dimensional faces form the structure maps of an operad \cite{jds}.  However, Stasheff and Markl show that the inclusion maps above give the collection $\{W_n\}$ a right module operad structure \cite[\S2]{mar}.
\end{rem}

\begin{figure} [h]
\centering {\includegraphics {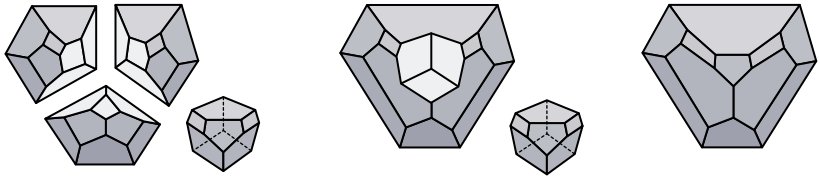}}
\caption{Four copies of $K_5$ in $W_4$}
\label{w4parts}
\end{figure}

\begin{rem}
The inclusion $K_n \hookrightarrow W_n$ coming from the proposition above shows associahedra as {\em faces} of the cyclohedron. A.\ Tonks has constructed an explicit map $W_n \rightarrow K_{n+1}$ between polytopes of the {\em same} dimension, in the spirit of his map from the permutohedron $P_n$ to $K_{n+1}$ \cite{ton}.  It is also noteworthy to point out a similar insight discovered by A.\ Ulyanov, where he shows associahedra gluing together to form the cyclohedron.

\begin{prop} \textup{\cite{uly}}
The cyclohedron \,$W_n$ is made up of \,$n$ copies of \,$K_{n+1}$.
\end{prop}

\noindent For example, the line segment $W_2$ is made up of two $K_3$ line segments with a pair of end points identified.  Figure~\ref{w4parts} shows the example of how $W_4$ can be constructed using four copies of $K_5$.
\end{rem}

\subsection{}
The original definition of $W_n$ describes it in terms of a compactified configuration space of {\em unlabeled} points \cite{bt}.  However, by adding labels, we go from one polytope to a space tiled by copies of it.

\begin{defn}
The moduli space \,\Z{n} is a real {\em Fulton-MacPherson}\footnote{Bott and Taubes use the Axelrod-Singer compactification \cite{as}, which is the differential geometric version of Fulton-MacPherson differing by a $\Int_2$ quotient at the blow-ups.} compactification of the configurations of $n$ {\em labeled} points on $S^1$ quotiented by the action of $S^1$. That is,
$$\Z{n} = \overline{{\rm Config}^n(S^1)/S^1}.$$
\end{defn}

\begin{figure}[h]
\centering {\includegraphics {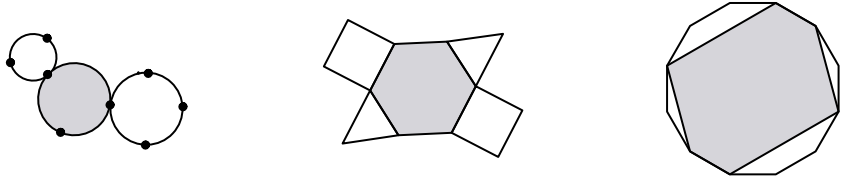}}
\caption{Bubbles and Centrally Symmetric Polygons}
\label{bubble2gon}
\end{figure}

The manifold \Z{n} has dimension $n-1$ and is without boundary.  It is a point when $n=1$ and $S^1$ when $n=2$.  As we will show, it becomes non-orientable for $n>2$.  Recall the difference in compactifications noted in \S\ref{ss:compact}.  The moduli space we have defined  using the Fulton-MacPherson method distinguishes the basespace $S^1$ from the bubbles $\RP$ coming from collisions.   We keep track of this information by shading in the circle corresponding to the base space, shown in  Figure~\ref{bubble2gon}, along with its duality in terms of $2n$-gons with centrally symmetric diagonals.  It naturally demonstrates the appearance of $W_n$ as its fundamental domain.  Note that whenever points collide, this is represented as drawing a diagonal {\em and} its symmetric counterpart.

One might ask why a $2n$-gon is needed in the case of \Z{n}.  The reason is due to the group action on the configuration space.  Since the symmetric group $\Sg_3$ is in $\PGL$, this allows us to fix three out of the $n$ points allowed to move on any of the bubbles in \M{n}.  In terms of polygon duals, we are guaranteed at least a triangle.  For \Z{n}, although there is still an action of $\PGL$ on the bubbles, the base space only has $S^1$ acting on it.  To help visualize the properties in terms of polygons, we double the sides of the base space (giving a $2n$-gon) and take care of the bubblings by adding centrally symmetric diagonals to the polygon.


\subsection{}
\begin{figure} [h]
\centering {\includegraphics {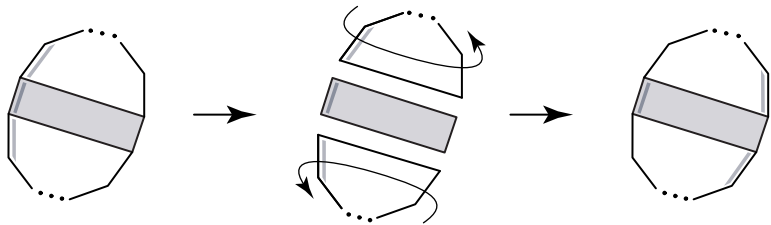}}
\caption{{\em Twist} along a diagonal}
\label{wtwist}
\end{figure}

It was the combinatorial construction of \M{n} using the twist operation \cite[\S3]{dev} which prompted Simion to wonder about a space of cyclohedra.  Although we use her Definition~\ref{d:cyclo} of the cyclohedron, our results are independent of her work. We mimic the idea above to show how cyclohedra glue together.  The {\em twist} of a $2n$-gon $G$ along a (pair of centrally symmetric) diagonal $d$, denoted by $\nabla_d(G)$, is the polygon obtained by separating $G$ along $d$, `twisting' (reflecting) both pieces symmetrically, and gluing them back (Figure~\ref{wtwist}).  Context will make it clear whether twist is being performed on an $n$-gon or a symmetric $2n$-gon.

We develop the combinatorial formation of \Z{n} as follows:  Label the edges of a $2n$-gon with $\{1, 2, \ldots, n\}$ such that the antipodal edge has the same label.  There will be $n!$ such possible labelings. However, since \Z{n} is defined modulo the action of $S^1$, two such labeled polygons are identified by a rotational action, yielding $(n-1)!$ distinct labeled polygons.  We obtain the following:
 
\begin{thm}
The space \,\Z{n} is tiled by $(n-1)!$ copies of \,$W_n$.  Two centrally symmetric labeled $2n$-gons \,$G_1, G_2$ with $k$ sets of diagonals (representing codim $k$ faces of cyclohedra) are identified in \,\Z{n} if there exist diagonals $d_1, \ldots, d_r$ of \,$G_1$ such that $$(\nabla_{d_1} \cdots \nabla_{d_r}) (G_1) = G_2.$$
\end{thm}

\begin{figure} [h]
\centering {\includegraphics {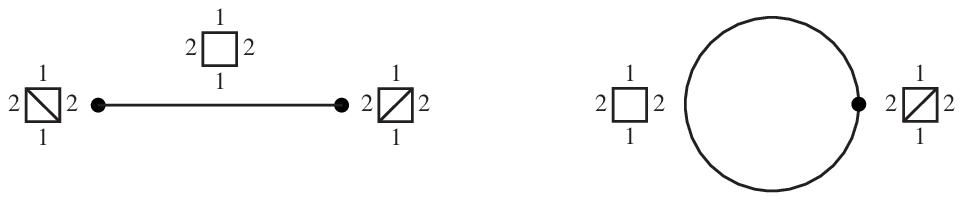}}
\caption{$W_2$ and \Z{2}}
\label{q02pcs}
\end{figure}

\noindent The proof of is an immediate consequence of the argument given in \cite[\S3.1]{dev}. Roughly, since $\PGL$ acts on each bubble and since $D_n \subset \PGL$, there is a reflection (twisting in terms of polygon duals) allowed enabling polytopes to be identified.  Figure~\ref{q02pcs} demonstrates the construction of \Z{2} from one copy of $W_2$ to form $S^1$.
Figure~\ref{q03pcs}b shows the gluing of two $W_3$ hexagons resulting in \Z{3}.

\begin{rem}
Faces of the associahedron are decomposed into the categories of dimension, type, and class, where they are used to better understand the structure of \M{n} (see \cite{dr} for definitions and details).  For example, the enumeration of faces with respect to dimension leads to the Euler characteristic of \M{n}.  Using methods and results from \cite{dr}, the initiated reader could easily carry out similar combinatorial calculations for the cyclohedron to further describe properties of \Z{n}.
\end{rem}

\begin{rem}
The classical $A_{\infty}$ operad is defined from the associahedron.  The map
$$K_n \times \Sg_n \rightarrow \dM{n+1},$$
where \dM{n+1} is the double cover with one point marked with a tangent direction, enables us to define a cyclic operad of {\em mosaics}, having \M{n} as its underlying space \cite[\S2]{dev}.  Although algebras over the $A_{\infty}$ operad are not commutative, those over the mosaic operad in some sense are, defining a class of {\em coherently homotopy associative commutative algebras}.
Similarly, the module operadic structure of $W_n$ transfers to \Z{n} with additional cyclic properties coming from the surjection
$$W_n \times \Sg_{n-1} \rightarrow \Z{n}$$
given by the twist operation.
\end{rem}

%
%
\section{Coxeter Groups and Blow-Ups}

\subsection{}
In light of understanding the global structure of these moduli spaces, the relationship between Coxeter groups and the polytopes $K_n$ and $W_n$ is introduced. We begin by looking at certain hyperplane arrangements and refer the reader to Brown \cite{bro} for any underlying terminology.

\begin{defn}
For a Euclidean space $V$, a {\em linear hyperplane} is a subspace $H \subset V$ of codim one passing through the origin.  A {\em finite reflection group} $W$ is a finite group of linear transformations of $V$ generated by reflections over a set of linear hyperplanes.
\end{defn}

\begin{figure} [h]
\centering {\includegraphics {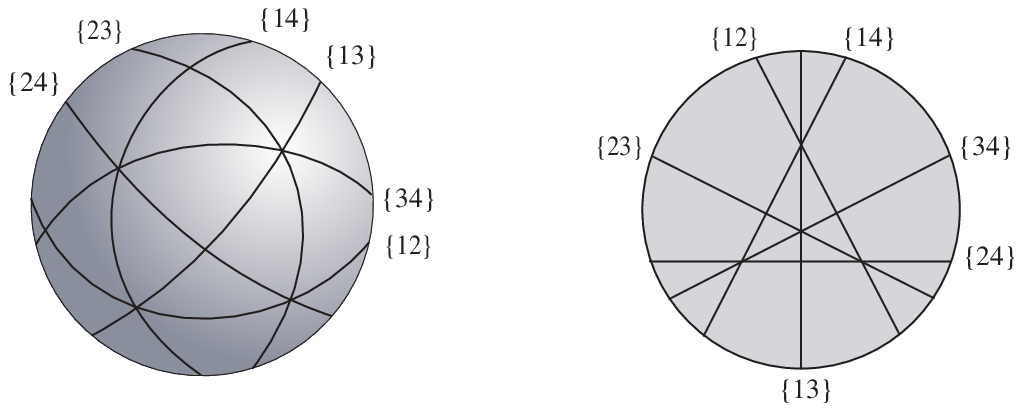}}
\caption{Braid arrangement on $\Sg V^3$ and $\Pj V^3$}
\label{svpv}
\end{figure}

Let $V^n \subset \R^{n+1}$ be the hyperplane defined by $\Sigma x_i = 0$ and $W = \Sg_{n+1}$.  Then $W$ is generated by transpositions acting as orthogonal reflections across the hyperplanes $x_i = x_j$.  The collection of such hyperplanes forms the {\em braid arrangement}.  Let $\Sg V^n$ and $\Pj V^n$ respectively be the sphere and the projective space in $V^n$; that is, $\Pj V^n$ is isomorphic to $\R\Pj^{n-1}$. Note that the braid arrangement gives $\Pj V^n$ a CW-cellular decomposition into $\frac{1}{2}(n+1)!$ open simplices. Figure~\ref{svpv} depicts the $n=3$ case.

\begin{defn}
For a linear subspace $X \subset Y$, we {\em blow-up \,$\Pj Y$ along \,$\Pj X$} by removing $\Pj X$, replacing it with the sphere bundle associated to the normal bundle of $\Pj X \subset \Pj Y$, and then projectifying the bundle.
\end{defn}

Blow-ups are the tools used by the compactification, keeping records of directions in which the colliding points approach $X$.  A general collection of blow-ups is usually non-commutative in nature; in other words, the order in which spaces are blown up is important.  For a given arrangement, De Concini and Procesi \cite{dp} establish the existence (and uniqueness) of a {\em minimal building set}, a collection of subspaces for which blow-ups commute for a given dimension, and for which the resulting manifold has only normal crossings.  For the braid arrangement, we use a combinatorially equivalent definition given in \cite[\S4.1]{dev}.

\begin{defn}
A  codim $k$ {\em minimal} element $\m_k$ of the minimal building set of $\Pj V^n$ is a subspace where $\binom{k+1}{2}$ hyperplanes of the braid arrangement intersect.
\end{defn}

\noindent It is due to M.\ Kapranov that the connection between blowing up the braid arrangement and the moduli space \M{n} is made.  An intuition for it comes from the definition of the braid arrangement, Theorem~\ref{t:twist}, and the construction of the associahedron using truncations.

\begin{thm} \textup{\cite[\S4.3]{kap}}
The iterated blow-up of \,$\Pj V^n$ along the cells $\{\m_k\}$ in \,{\em increasing} order of dimension yields \,\M{n+2}.
\end{thm}

\begin{figure} [h]
\centering {\includegraphics {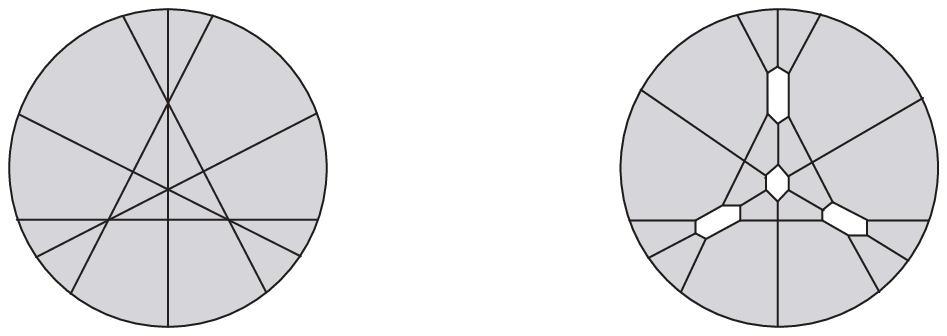}}
\caption{$\Pj V^3$ and \M{5}}
\label{pvm05}
\end{figure}

\begin{exmp}
Figure~\ref{pvm05} shows $\Pj V^3$ before and after blowing up.   The minimal elements of codim two turn out to be the points of triple intersection; blowing up along these components yield hexagons with antipodal identifications.  The diagram on the left has $\R\Pj^2$ tiled by twelve simplices; the one on the right after blow-ups is the connected sum of five projective planes, tiled by twelve $K_4$ pentagons.
\end{exmp}

\begin{rem}
The diagram shown in Figure~\ref{pvm05}b is first found in the work of Brahana and Coble in $1926$~\cite[\S1]{bc}, relating to possibilities of maps with twelve five-sided countries.
\end{rem}


\subsection{}
As demonstrated, blowing up minimal elements of arrangements has a local consequence of truncating the simplex to construct the associahedron.  This is best understood in light of Coxeter groups.  The Coxeter (or Weyl) group of type $A_n$ is the symmetric group $\Sg_{n+1}$; moreover, it is the full symmetry group of the $(n+1)$-simplex, coming from representing the simplex as a convex linear combination of vectors.  The Coxeter diagram of $A_n$
$$\bullet - \bullet \cdots \bullet - \;\bullet$$
with $n$ nodes gives a presentation by $n$ transpositions $\{s_i\}$, with the relations  $(s_i,s_{i+1})^3 = 1$ for adjacent generators and $(s_i, s_j)^2=1$ otherwise.  From the braid arrangement, we can identify the nodes of the Coxeter diagram of $A_n$ with the top dim faces of the $(n-1)$-simplex. 

\begin{prop}
Truncating the faces of the simplex which correspond to connected subdiagrams of the Coxeter diagram of \,$A_n$ in \,{\em increasing} order of dimension results in \,$K_{n+1}$.
\end{prop}

\begin{figure} [h]
\centering {\includegraphics {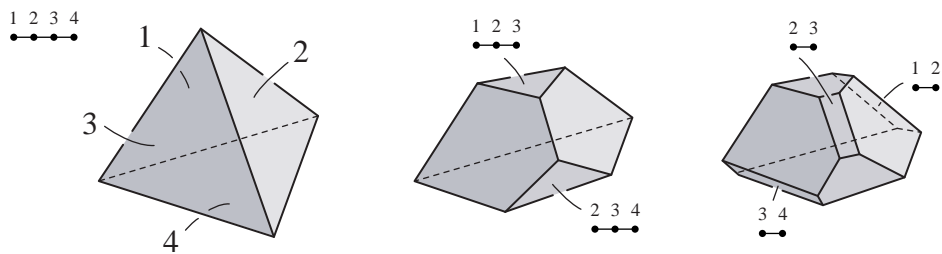}}
\caption{Truncation of simplex to $K_5$}
\label{trunk-k5}
\end{figure}

\noindent The proof follows immediately from the construction of $K_n$ given in \cite[\S5.4]{dev}.  Lee~\cite[\S3]{lee} also obtains this using a sequence of stellar subdivisions on the polar dual to the associahedron.  Figure~\ref{trunk-k5} shows an example of the $3$-simplex, with the initial truncation of vertices and then of edges, resulting in $K_5$.

\subsection{}
We now show the moduli space \Z{n} as an affine analogue of \M{n}.   Let $\widetilde A_n$ denote the {\em affine} root system of type $A_n$ (see Bourbaki \cite{bou} for details).   The Coxeter diagram of $\widetilde A_n$  is $n$ vertices arranged on a circle.  Its Coxeter group is the semi-direct product
$$\widetilde \Sg_n \,=\, \Int^{n-1} \,\rtimes\, \Sg_n\, ,$$
where the $n$ generators (transpositions) satisfy $(s_i,s_j)^3 = 1$ when adjacent on the diagram, and commute otherwise.  Let the nodes on the Coxeter diagram of $\widetilde A_n$ correspond to the top dimensional faces of the $(n-1)$-simplex. 

\begin{prop} \label{p:trunkw}
Truncating the faces of the simplex which correspond to connected subdiagrams of \,$\widetilde A_n$ in \,{\em increasing} order of dimension results in \,$W_n$.
\end{prop}

\begin{figure} [h]
\centering {\includegraphics {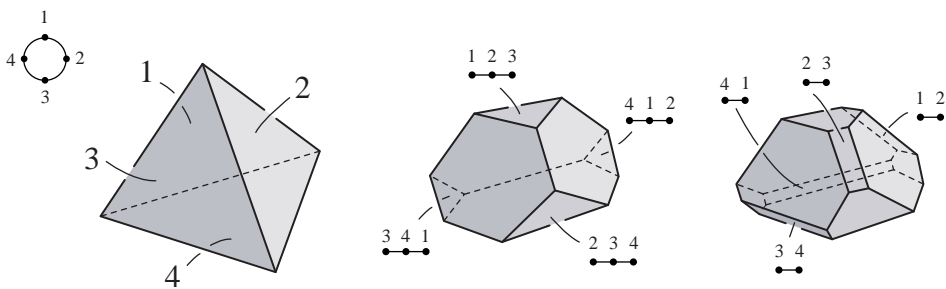}}
\caption{Truncation of simplex to $W_4$}
\label{trunk-w4}
\end{figure}

\noindent The proof of the statement follows from the construction of $W_n$ given in  \cite[Appendix B]{jds2}.  Figure~\ref{trunk-w4} shows an example from $\widetilde A_4$, with the initial truncation of vertices and then of edges, resulting in $W_4$.

In terms of hyperplane arrangements, we move from {\em spherical} geometry coming from linear hyperplanes to {\em Euclidean} geometry arising from affine hyperplanes.

\begin{defn}
An {\em affine hyperplane} of $V$ is a subset of the form $x + H$, where $x \in V$ and $H$ is a linear hyperplane of $V$. An {\em affine reflection group} $W$ is a locally finite group of affine isometries of $V$ generated by reflections over a set of affine hyperplanes.
\end{defn}

\begin{figure} [h]
\centering {\includegraphics {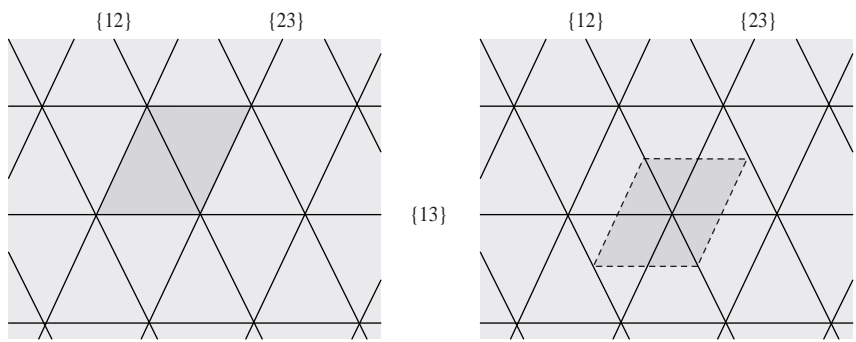}}
\caption{Tiling of $\R^2$ by hyperplanes}
\label{hextiles}
\end{figure}

We follow the construction given in Brown \cite[\S6.1]{bro}. Let $V^n \subset \R^n$ be the hyperplane defined by $\Sigma x_i = 0$.  The {\em affine} braid arrangement is the set of affine hyperplanes in $V^n$ of the form $x_i = x_j + k$, where $k \in \Int$, having $\widetilde \Sg_n$ as its affine reflection group.  This arrangement gives $V^n$ a CW-cellular decomposition into an infinite collection of open $(n-1)$-simplices. Figure~\ref{hextiles} gives the example when $n=3$ generated by the reflections $\{(12), (23), (13)\}$.

Previously, we looked to the projective sphere $\Pj V^n$ to give a compact space. In the affine case, the natural candidate is the quotient $\La V^n = V^n/\Int^{n-1}$, which is homeomorphic to the $(n-1)$-torus tiled by $(n-1)!$ simplices.  The shaded region in Figure~\ref{hextiles}, after quotienting, yields the two-torus $\La V^3$.  By changing the fundamental domain (Figure~\ref{hextiles}b), we see the braid arrangement $A_2$ naturally appearing in $\widetilde A_3$.  The affine braid arrangement also has a combinatorial notion of minimal building set: A codim $k$ minimal element \,$\m^k$ of the minimal building set of \,$\La V^n$ is a subspace where $\binom{k+1}{2}$ hyperplanes intersect. 

\begin{thm}
The iterated blow-up of \,$\La V^n$ along the cells \,$\{\m^k\}$ in \,{\em increasing} order of dimension yields \,\Z{n}.
\end{thm}

\begin{proof}
We show the moduli space \Z{n} as the minimal blow-up of the Coxeter group $\widetilde A_{n}$ modulo the action of $\Int^{n-1}$ inside the translation subgroup.   Note that there are $(n-1)!$ simplices in $\La V^n$, each with a particular ordering given by the hyperplanes.  It follows from Proposition~\ref{p:trunkw} that the iterated blow-ups of $\La V^n$ do indeed give a tiling by cyclohedra $W_n$.  The affine braid arrangement is simply a way of keeping track of all the gluings defined by the twisting operation in a global setting.  Adjacent polytopes can be accessed by crossing a hyperplane ($x_i = x_j$) or an intersection of hyperplanes, thereby permuting their labelings.  Indeed, our notation of labeled polygons with twisting of diagonals mimics the information found in the hyperplane arrangement. 
\end{proof}

\begin{figure} [h]
\centering {\includegraphics {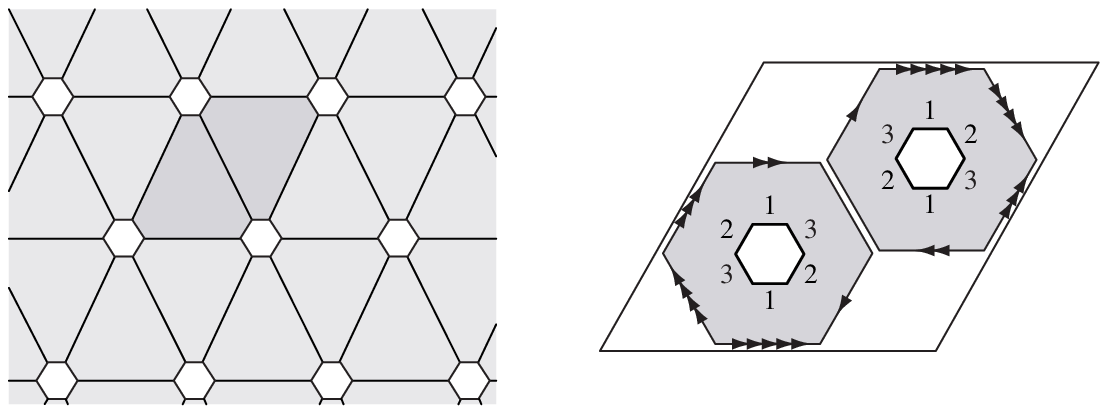}}
\caption{Copies of $W_3$ gluing to form \Z{3}}
\label{q03pcs}
\end{figure}

\begin{exmp}
Figure~\ref{q03pcs}a shows the result after blowing up the minimal elements of $\La V^3$.  We see a tiling by the hexagons $W_3$, illustrated in detail on the right of the figure.  The resulting manifold \Z{3} is homeomorphic to $\R\Pj^2 \# \R\Pj^2 \# \R\Pj^2$.
\end{exmp}

\subsection{}
Our motivation in further analyzing \Z{n} comes from some observations in Figure~\ref{q03pcs}.  The minimal point after blow-up becomes a hexagon with antipodal identification.  This is essentially a triangle, which can be seen as \M{4} shown in Figure~\ref{m04pcs}b.  Note also that the three lines on $\La V^3$ become line segments with endpoint identifications in \Z{3}, which is seen from Figure~\ref{q02pcs} as three copies of \Z{2}. In general, we observe the following:

\begin{thm}
In \ \Z{n}, there are \ $\binom{n}{k+1}$\, copies of \ $\M{k+2} \times \Z{n-k}$.
\end{thm}

\begin{proof}
We show that the blow-up of a minimal element $\m_k \in \La V^n$ results in $\M{k+2} \times \Z{n-k}$.  From the combinatorics of the affine braid arrangement, it is not hard to see that each $\m_k$ corresponds to a choice of $k+1$ elements from the set $\{1, \ldots, n\}$.  Choose an arbitrary minimal cell and assign it such a choice, say $\{p_1, \ldots, p_{k+1}\}$.  We view this as a centrally symmetric $2n$-gon having a diagonal partitioning it into a pair of $k+1$ symmetrically labeled sides $\{p_1, \ldots, p_{k+1}\}$, with the sides of the central polygon using the remaining labels.  Note the correspondence of this dissected $2n$-gon to the product $K_{k+1} \times W_{n-k}$.  As this minimal cell is blown-up, we see $(k+1)!$ different ways in which the labels can be arranged on the non-central polygons.  However, since {\em twisting} is allowed along the diagonals, we get $\frac{1}{2}(k+1)!$ different labelings, each corresponding to a $K_{k+1}$.  Indeed, this is {\em exactly} how one gets \M{k+2}, with the associahedra glued as defined above.  Therefore, a fixed labeling of the central polygon gives $\M{k+2} \times W_{n-k}$, while allowing all possible labelings gives the result.
\end{proof}

\begin{figure} [h]
\centering {\includegraphics {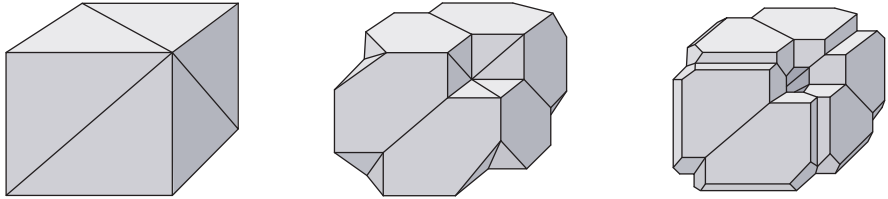}}
\caption{Blow-ups leading to six $W_4$ polytopes tiling \Z{4}}
\label{q04pcs}
\end{figure}

\begin{exmp}
We move to the $n=4$ case to further illustrate the theorem above. The left most drawing of Figure~\ref{q04pcs} shows the fundamental cube tiling $\R^3$ with hyperplane cuts. Blowing up the minimal point (codim three) and then the codim two elements yield the middle and the right hand pictures.  The space \Z{4} is shown to be tiled by $3!$ copies of $W_4$.  Notice the vertex of the cube (identified to form the $3$-torus) becomes \M{5} after blow-ups, tiled by twelve pentagons.  The minimal lines become triangular tori, that is, $\M{4} \times \Z{2}$.  The faces of the cube can be seen to form the $2$-torus with blow-ups, resulting in \Z{3} tiled by two $W_3$ hexagons.
\end{exmp}

%
%
\section{Some Observations}

\subsection{}
We have shown the cyclohedron to be the $\widetilde A_n$ affine analogue of the associahedron. However, Simion thought of it as a type $B_n$ associahedron.\footnote{By taking minimal blow-ups coming from the hyperoctahedral root system of type $B_n$, we would still arrive at the associahedron, not the cyclohedron.}  The reasoning behind this is explained by looking at the combinatorial question that motivated her.  We thank V.\ Reiner for helpful insights into his work \cite{rei}.

\begin{defn}
A {\em $k$-partition} separates $\{1, \ldots, n\}$ into disjoint sets $B_1, \ldots, B_k$.  A partition is {\em non-crossing} if for elements $x_1, x_2 \in B_i$ and $y_1, y_2 \in B_j$ such that $x_1 < y_1 < x_2 < y_2$, then $i=j$.  Let $NC(n,k)$ be the number of non-crossing $k$-partitions of $\{1, \ldots, n\}$.
\end{defn}

\noindent The $f$-vector $\{f_i\}$ of a polytope $P$ is a combinatorial invariant, where $f_k(P)$ equals the number of $k$-dimensional faces.  Another method of encapsulating this information is by the $h$-vector $\{h_i\}$ (see \cite{zie} for details). Denoting $K_n^*$ as the combinatorial dual of the associahedron, the following remarkable equivalence exists:
$$NC(n,n-k) = h_k(K_{n+1}^*).$$

The poset of partitions of $\{1, \ldots, n\}$ is isomorphic to the lattice of intersection subspaces for the type $A_n$ braid arrangement.  For example, a partition $\{1,2,4\}\{3,5\}$ corresponds to the intersection subspace of $x_1=x_2=x_4$ and $x_3=x_5$. Similarly, by replacing the set with $\{1, \ldots, n, \bar 1, \ldots, \bar n\}$, we can consider the hyperoctahedral hyperplane arrangement for type $B_n$.  The barring is an involution denoting the sign change coming from the hyperplanes $x_i = \pm x_j$.

\begin{defn}
A {\em $k_B$-partition} separates $\{1, \ldots, n, \bar 1, \ldots, \bar n\}$ into disjoint sets $B_0, B_1, \ldots, B_k, \bar B_1, \ldots, \bar B_k$, where $\bar B_i = \{\bar x \ |\ x \in B_i \}$ and $B_i = \bar B_i$ only when $i=0$.
\end{defn}

\noindent  Let $NC_B(n,k)$ be the number of non-crossing $k_B$-partitions of $\{1, \ldots, n, \bar 1, \ldots, \bar n\}$. Simion was able to show that the cyclohedron satisfies the equality:
$$NC_B(n,n-k) = h_k(W_{n+1}^*).$$

\begin{rem}
Combinatorially, Simion was able to construct a space tiled by eight copies of $W_3$ hexagons, using the twist operation defined in \cite{dev}.  We now see that she had arrived at the four-fold cover of \Z{3} (in Figure~\ref{hextiles}, choose a region four times that of the shaded one).  By labeling the antipodal sides of the polygons with $i$ and $\bar i$, a $2^n$-fold cover of the moduli space \Z{n} is created.  
\end{rem}

\begin{rem}
Burgiel and Reiner \cite{br} discuss two more type $B_n$ analogues of the associahedron, both being different from the cyclohedron.
\end{rem}


\subsection{}
The cyclohedron made its debut in knot theory implicitly through the work on Kontsevich \cite{kon} and explicitly by Bott-Taubes.  The imbedding $S^1 \rightarrow \R^3$ defining a knot induces a map onto their compactified $n$ point configuration spaces, thereby introducing $W_n$.  Roughly, one  then `carefully' integrates certain pullbacks of volume forms on $S^2$ using the inclusion $S^1 \rightarrow \R^3 \rightarrow S^2$; for a good topological understanding, see the work of D.\ Thurston~\cite{thu}.

The cyclohedron also appears in the work of Bar-Natan in a different setting \cite[\S4]{bar}.  He gives a combinatorial approach of sketching some relations that arise in computing certain Kontsevich-KZ invariants coming from braids and chord diagrams (note that chord diagrams come from $2n$ marked points on $S^1$ with pairwise identification).  He introduces a multitude of relations and their corresponding polytopal structures, including Kapranov's permutoassociahedron $KP_n$ and $W_n$. 


\subsection{}
Davis, Januszkiewicz, and Scott show the minimal blow-ups of certain hyperplane arrangements to be $K(\pi, 1)$ spaces \cite[\S5]{djs}.  In other words, the homotopy properties of these spaces are completely encapsulated in their fundamental groups. The crucial result motivating their work comes from Gromov \cite{grm}, where he relates nonpositive curvature of manifolds to flag complexes.  For our purposes, we get the following striking fact:

\begin{thm}
The moduli spaces \,\M{n} and \,\Z{n} are aspherical.
\end{thm}

\begin{figure} [h]
\centering {\includegraphics {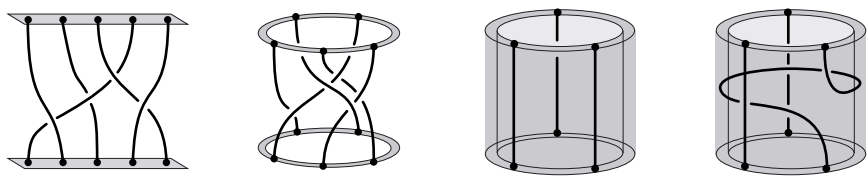}}
\caption{Four braids, two coming from $A_4$ and two from $\widetilde A_2$}
\label{braids}
\end{figure}

There is a close similarity between the fundamental group of \M{n} and the braid group, described from a combinatorial viewpoint in \cite[\S6]{dev}.  Recall that the {\em Artin group} of $A_n$ is the classic braid group on $n+1$ strings between points on an interval (Figure~\ref{braids}a) or on a circle (Figure~\ref{braids}b); simply cut the circles open and lay the strings flat to obtain an isomorphism.  The affine analog of classical braids is the Artin group of $\widetilde A_n$; this can be presented as braids on $n+1$ strands restricted to a cylindrical shell. The latter diagrams in Figure~\ref{braids} show two distinct elements of the $\widetilde A_2$ braid group. There is a parallel between these cylinder braids and the fundamental group of \Z{n}, with notably deeper relations to knot invariants and the cyclohedron.

\subsection{}
It seems quite natural to generalize what we have done.  The motivation for the real moduli space of points \M{n} comes from its {\em complex} counterpart \CM{n}.  The importance of \CM{n} can be seen, for example, as being a fundamental building block in the theory of Gromov-Witten invariants, among numerous other appearances in literature.  A reasonable candidate for the complex analogue for \Z{n} could be:
$$\Z{n}(\C) = \overline{{\rm Config}^n(S^2)/SO(3)}.$$
Again, we would like to know the role this space would play in mathematical physics and algebraic geometry. Similarly, we can generalize by mimicking the construction to other Coxeter groups.  Much has been done in terms of the minimal blow-ups in \cite{djs}.  However, there have been a wealth of new compactifications that have been studied recently and we discuss the interplay between these ideas in a forth coming work \cite{dev2}.

%
%

\bibliographystyle{amsplain}

\end{document}